\documentclass[cmp,referee]{svjour}
\usepackage{amsmath}
\usepackage{amssymb}
\makeatletter
\usepackage{amsmath}
\usepackage{amsfonts}
\usepackage[normalem]{ulem}
\usepackage{txfonts}
\usepackage{theorem}
\usepackage{makeidx}
\usepackage{amsmath}
\usepackage{epic}
\usepackage{amscd}
\usepackage{bbm}
\usepackage{xy}
\usepackage{times}
\usepackage{verbatim}
\usepackage{graphics,epsfig}
\usepackage{bm}

\@addtoreset{equation}{section}
\makeatother
\newcommand{\ue}{\mathrm{e}}
\newcommand{\ui}{\mathrm{i}\,}

\newcommand{\aot}{A^{\scriptscriptstyle 1/2}}
\newcommand{\kd}{{k_{\rm d}}}
\newcommand{\kt}{{\tilde k}}
\newcommand{\sfa}{{{\sf a}(\kappa)}}

\newcommand{\sfd}{{{\sf d}(\kappa)}}
\newcommand{\sfe}{{{\sf e}(\kappa)}}
\newcommand{\sff}{{{\sf f}(\kappa)}}

\begin{document}


\title{Entire solutions of hydrodynamical equations with exponential 
dissipation}
\author{Claude Bardos\inst{1}, Uriel Frisch\inst{2}, Walter Pauls\inst{3}, 
Samriddhi Sankar Ray\inst{4}, and Edriss S. Titi\inst{5} } 
\institute{Universit\'e Denis Diderot and Laboratoire J.L. Lions Universit\'e Pierre et Marie Curie, 
Paris, France
\and
UNS, CNRS, Laboratoire Cassiop\'ee, OCA, 
BP 4229, 06304 Nice cedex 4, France
\and
Max Planck Institute for Dynamics and Self-Organization, G\"ottingen, Germany
\and
Center for Condensed Matter Theory, 
Department of Physics, Indian Institute of Science, 
Bangalore, India 
\and 
Department of Mathematics and
Department of Mechanical
and Aerospace Engineering, 
University of Irvine,
CA 92697-3875, USA
\\
Department of Computer Science 
and Applied Mathematics, 
Weizmann Institute of Science, 
Rehovot 76100, Israel
}
\communicated{name}
\date{Received: date / Accepted: date}
\maketitle
\begin{abstract}
We consider a modification of the three-dimensional Navier--Stokes
equations and other hydrodynamical evolution equations with
space-periodic initial conditions in which the usual Laplacian of the
dissipation operator is replaced by an operator whose Fourier symbol
grows exponentially as $\ue ^{|k|/\kd}$ at high wavenumbers $|k|$.
Using estimates in suitable classes of analytic functions, we show
that the solutions with initially finite energy become immediately
entire in the space variables and that the Fourier coefficients decay
faster than $\ue ^{-C(k/\kd) \ln (|k|/\kd)}$ for any $C<1/(2\ln
2)$. The same result holds for the one-dimensional Burgers equation
with exponential dissipation but can be improved: heuristic arguments
and very precise simulations, analyzed by the method of asymptotic
extrapolation of van der Hoeven, indicate that the leading-order
asymptotics is precisely of the above form with $C= C_\star
=1/\ln2$. The same behavior with a universal constant $C_\star$ is
conjectured for the Navier--Stokes equations with exponential
dissipation in any space dimension. This universality prevents the
strong growth of intermittency in the far dissipation range which is
obtained for ordinary Navier--Stokes turbulence.  Possible
applications to improved spectral simulations are briefly discussed.

\end{abstract}

\section{Introduction}
\label{s:intro}

More than a quarter of a millenium after the introduction by Leonhard
Euler of the equations of incompressible fluid dynamics the question
of their well-posedness in three dimensions (3D) with sufficiently
smooth initial data is still moot
\cite{majda-bertozzi,fmb03,bardostiti,constantin07} (see also many papers
in \cite{efms} and references therein). Even more vexing is the fact
that switching to viscous flow for the solution of the Navier--Stokes
equations (NSE) barely improves the situation in 3D \cite{lions,constantinfoias,fefferman,temam,sohr}.
Finite-time blow up of the solution to the NSE can thus not be ruled
out, but there is no numerical evidence that this happens.

In contrast, there is strong numerical evidence that for analytic
spatially periodic initial data both the 3D Euler and NSE
have complex space singularities. Indeed, when such equations are
solved by (pseudo-)spectral techniques the Fourier transforms of the
solution display an exponential decrease at high wavenumbers,  which
is a signature of complex singularities \cite{brachetetal}. This behavior was
already conjectured by von Neumann \cite{vonneumann}  
who pointed out on p.~461 that
the solution should be  analytic with an exponentially
decreasing spectrum.  Recently Li and Sinai used a
Renormalization Group method to prove that for certain complex-valued
initial data the 3D NSE display finite-time blow up in
the real domain (and, as a trivial corollary, also in the complex
domain) \cite{LS07}. 

For some PDEs in lower space dimensions explicit
information
about the position and type of complex singularities may be available.
For example, complex singularities can sometimes be related to poles of elliptic
functions in connection with  the reaction diffusion equation
\cite{olivertiti} and  2D incompressible Euler
equations 
in Lagrangian coordinates \cite{paulsmatsumoto}. The best understood case
is that of the  1D Burgers equation with ordinary (Laplacian) dissipation:\footnote{The case of the Burgers equation
with modified dissipation will be considered in
Section~\ref{s:burgers}.} its singularities are poles located at the zeroes
of the solutions
of the heat equation to which it can be mapped by using the Hopf--Cole
transformation (see, e.g., \cite{senoufetal,polaciksverak} and references
therein).

We now return to the 3D NSE with real analytic data. It is known
that blow up in the real domain can be avoided altogether by modifying the dissipative operator, whose Fourier-space symbol
is $\mu |k|^2$, to a higher power of the Laplacian with symbol $\mu
|k|^{2\alpha}$ ($\alpha >5/4$) \cite{lions,ladyzhenskaya}. The numerical evidence is however that
complex singularities cannot be avoided by this ``hyperviscous''
procedure, frequently used in geophysical simulations (see, for
example, \cite{holloway}).

Actually, we are unaware of any instance of a nonlinear space-time
PDE, with the property that the Cauchy problem is well posed in the
complex space domain for at least some time and which is guaranteed
never to have any complex-space singularities at a finite distance
from the real domain. In other words the solution 
stays or becomes entire for all $t>0$. Here we shall
show that solutions of the Cauchy problem are entire for a fairly
large class of pseudo-differential nonlinear equations, encompassing variants of the 3D NSE, which possess ``exponential dissipation'', that is
dissipation with a symbol growing exponentially as $\ue^{|k|/\kd}$ 
with the ratio of the wavenumber $|k|$ to a reference wavenumber $\kd$.

The paper is organized as follows. In Section~\ref{s:entire} we
consider the forced 3D incompressible NSE  in a periodic
domain with exponential dissipation. The initial conditions are
assumed  just to have finite energy. The main theorem is
established using classes of analytic functions whose norms contain
exponentially growing weights in the Fourier space 
\cite{foiastemam,ferrarititi}.  In
Section~\ref{s:decay} we show that the Fourier transform of the
solution decays at high wavenumbers faster than $\exp{(-C\kt\ln \kt)}$ for any
$C< 1/(2\ln 2)$. Here, $\kt \coloneqq |k|/\kd$ is the nondimensionalised
wavenumber.
In Section~\ref{s:extensions} we briefly present
extensions of the result to other instances: different space dimensions and 
dissipation rates, problems formulated in the whole space and on a sphere, and different equations. In Section~\ref{s:burgers} we then 
turn to the 1D
Burgers equation with a  dissipation growing exponentially at high 
wavenumbers, for which the same bounds
hold as for the 3D Navier--Stokes case. However in the
Burgers case, simple heuristic considerations (Section~\ref{ss:heuristics}) and very accurate numerical simulations performed by two
different techniques (Sections~\ref{ss:simulation1} and \ref{ss:simulation2}), indicate that the leading-order asymptotic decay is
precisely $\exp{((-1/\ln 2)\kt\ln \kt)}$. We observe that the heuristic
approach, which involves a dominant balance argument applied 
in spatial Fourier space,
is also applicable to the  3D Navier--Stokes case with exactly the
same prediction regarding the asymptotic decay. In the concluding
Section~\ref{s:conclusion} we discuss open problems and a possible
application.

\section{Proof that the solution is entire}
\label{s:entire}

We consider the following 3D spatially periodic Navier--Stokes equations with an exponential
dissipation (expNSE)
\begin{eqnarray}
&&\frac{\partial u}{\partial t} + u\cdot\nabla u = -\nabla p - \mu {\cal D} u +
f, \qquad \nabla\cdot u=0,\label{expnsoldfashion}\\
&&u(x,0)=u_0(x). \label{initoldfashion}
\end{eqnarray}
Here, ${\cal D}$ is the (pseudo-differential) operator whose
Fourier space symbol is $\ue^{2\sigma |k|}$, that is a dissipation rate
varying
exponentially with the wavenumber $|k|$, $u_0$ is the initial
condition, $f$ is a prescribed driving force and $\mu$ and $\sigma$
are prescribed positive coefficients. The problem is
formulated in a periodic domain $\Omega$ (for simplicity 
of notation we take $\Omega=[0,2\pi]^3$). The driving force is 
assumed to be a divergence-free trigonometric polynomial in the
spatial coordinates.
For technical convenience we use $\sigma$ in the statements
and
proofs of mathematical results, while the use of the reference
wavenumber $\kd = 1/(2\sigma)$ is preferred when discussing the
results.

The initial condition is  taken to be a divergence-free periodic vector
field
with a finite $L^2$ norm (finite
energy).

As usual the problem is rewritten as an abstract ordinary differential
equation in a suitable function space, namely 
\begin{eqnarray}
&&\frac{du}{dt} + \mu \ue^{2\sigma {\aot}} u+
B(u,u)= f; \label{expns}\\
&&u(0)=u_0;  \label{init}
\end{eqnarray}
where $A \coloneqq -\nabla ^2$ and $B(u,u)$ is a suitable quadratic form
which takes into account the nonlinear term, the pressure term and the
incompressibility constraint (see, e.g. \cite{lions,constantinfoias,temam}). Note that the Fourier symbol of
$A^{1/2}$ is $|k|$. 

The problem is formulated in the space
$H\coloneqq\{\varphi \in
(L^2(\Omega))^3:\varphi$ is periodic, $\int\varphi\,dx=0$, 
\hbox{$\nabla\cdot\varphi=0\}$.} Here, for any $\lambda\ge 0$,  the Fourier 
symbol of the operator $\ue^{\lambda
{\aot}}$is given by $\ue^{\lambda
|k|}$, where $k \in \mathbb{Z}^3 \backslash \{(0,0,0)\}$.

To prove the entire character, with respect to the spatial variables, of the solution $u(t)$ of expNSE for
$t>0$, it suffices to
show that its Fourier coefficients decrease faster than exponentially
with the wavenumber $|k|$. This will be done by showing that, for any
$\lambda >0$,  the $L^2$ norm of $\ue^{\lambda \aot}u$,   the solution
with an exponential weight in Fourier space, is finite. As usual, we here denote the
$L^2$ norm of  a real space-periodic function $f$ by $|f|\coloneqq \sqrt{
\int_{[0,2\pi]^3} |f(x)|^2 dx}$. Moreover, $H^m$ will be the usual
$L^2$ Sobolev space of index $m$ (i.e., functions which have up to $m$
space derivatives in $L^2$).\\

The main result (Theorem~2.1) will make use of the following
Proposition
which was inspired by \cite{foiastemam} (see also \cite{ferrarititi})\\
\textbf{Proposition 2.1} \,\,
{\it Let $\alpha \geq 0$  ,
$\beta > 0$, $\kappa > 3$
and $\varphi \in {\rm dom}\,(\ue^{(\alpha+\beta){\aot}})$. Then

\begin{equation}
\left |\ue^{\alpha {\aot}}B (\varphi, \varphi)\right |\:\ \leq
C_{\rm A} \left(l_\kappa(\beta)\right) ^\sfa\:\left |\ue^{\alpha
{\aot}}\varphi\right |^{2-\sfa}~~\left |\ue^{(\alpha+\beta)
{\aot}} \varphi\right |^\sfa,
\label{prop1}
\end{equation}
where $C_{\rm A}$ is a universal constant and
\begin{equation}
l_\kappa(\beta)\coloneqq \sup_{0\leq x<\infty} x^\kappa \ue^{-\beta x} =\left(
\frac{\kappa}{\beta} \right)^\kappa \ue^{-\kappa}.
\label{deflkappa}
\end{equation}\\
}
\textit{Notation}\,\, In Proposition~2.1 and also in the
sequel
we use 
the following notation (to avoid fractions in exponents):
\begin{eqnarray}
&& \sfa \coloneqq \frac{5}{2\kappa},\quad    \sfd \coloneqq \frac{2\kappa+5}{2\kappa-5} \nonumber \\
&&\sfe \coloneqq \frac{4\kappa-5}{2\kappa-5},\quad  \sff \coloneqq
  \frac{2\kappa}{2\kappa -5}.
\label{defsfag}
\end{eqnarray}

\textit{Proof} \,\,Let $w \in H$. By using the Fourier
representations $\varphi(x) = \sum_l  \ue ^{\ui l\cdot x}\hat\varphi_l$
and  $w(x) = \sum_k  \ue ^{\ui k\cdot x}\hat w_k$ and Parseval's
theorem, we have
\begin{equation}
\frac{1}{(2\pi )^3 } \, \left (\ue^{\alpha {\aot}}B (\varphi,
\varphi),w
\right )=\sum_{k \in
\mathbb{Z}^3\: k\neq0}~~\ue^{\alpha|k|}  \left(\sum _{l+m=k;
\:l,m\neq0}(\hat \varphi_{l}\cdot \ui m)\hat \varphi_{m}\right)\cdot \hat
w_k^* ,
\end{equation}
where the $*$ means complex conjugation.

Since $\ue^{\alpha|k|} \leq \ue^{\alpha|m|+\alpha|l|}$, when $k=l+m$,
we can estimate the absolute value of  the right-hand side from above as
\begin{equation}
\begin{split}
&
\leq \sum_{k \neq 0 }  \sum_{l+m=k, \ l,m \neq 0}
e^{ \alpha\vert l \vert  }  
\vert \hat{\phi } _l \vert \, \vert m \vert  e^{ \alpha\vert m \vert  }   
\vert  \hat{\phi } _m  \vert \, 
\vert \hat{w} _k  \vert   
\\
&
=\frac{1}{(2\pi )^3 } 
\int  \Phi (x) \Psi (x) W (x) dx  \leq  
\frac{1}{(2\pi )^3 } \Vert \Psi \Vert _{L^{\infty } }  \left | \Phi \right |  
\left | W \right | ,
\end{split}
\end{equation}
where the functions $\Phi (x) $, $\Psi (x) $ and $W(x)$ are given by
\begin{equation}
\Phi(x)=\sum_{l\neq0}\ue^{\alpha|l|}|\hat \varphi_l | \ue^{\ui  l\cdot x},
\label{phifourier}
\end{equation}
\begin{equation}
\Psi(x)=\sum_{m\neq0}|m| \ue^{\alpha|m|}|\hat \varphi_m| \ue^{\ui  m\cdot
x},
\end{equation}
and
\begin{equation}
W(x)=\sum_{k\neq 0}|\hat w_k|\ue^{\ui  k \cdot x},
\end{equation}
and the last inequality follows from the Cauchy--Schwarz inequality.

By Agmon's inequality \cite{agmon} (see also \cite{constantinfoias}) in 3D we have
\begin{equation}
\|\Psi\|_{L^{\infty}}\leq
C_{\rm A}\|\Psi\|_{H^1}^{\frac{1}{2}}\|\Psi\|_{H^2}^{\frac{1}{2}}\leq 
C_{\rm A}\left
|{\aot} \Psi\right |^{\frac{1}{2}}\left |A \Psi\right |^{\frac{1}{2}}=
C_{\rm A}\left | A \ue ^{\alpha\aot} \varphi\right|\,\left| A ^ {3/2}
\ue ^{\alpha \aot} \varphi \right | ,
\label{agmon}
\end{equation}
where $C_{\rm A}>0$ is a universal constant.
By using \eqref{phifourier}, \eqref{agmon} and the fact that $|W|= |w|$, we obtain
\begin{equation}
\left |(\ue^{\alpha {\aot}}B(\varphi, \varphi),w)\right |\leq C_{\rm A} \left |\ue
^{\alpha {\aot}}\varphi\right |\:\left |A\ue^{\alpha
  {\aot}}\varphi\right |^{\frac{1}{2}}\:\left
|A^{\frac{3}{2}}\ue^{\alpha {\aot}}\varphi\right |^{\frac{1}{2}}\:|w|.
\end{equation}
And by using the interpolation inequality between $L^2$ and
$H^{\kappa }$,
where $\kappa > 3$, we obtain\footnote{The simplest formulation is obtained for $\kappa
  =5$ but the optimization of the bound for the law of  decay in
  Section~\ref{s:decay} requires using arbitrary $\kappa$.}
\begin{equation}
\left |(\ue^{\alpha {\aot}}B(\varphi, \varphi),w)\right |\leq C_{\rm A} \left |\ue^{\alpha
 {\aot}}\varphi\right |^{   {2-\sfa}   }\:\left
 |A^{\frac{\kappa}{2}}\ue^{\alpha {\aot}}\varphi\right |^{  \sfa   }\:|w|.
\end{equation}
Now, to obtain the inequality in Proposition~2.1, we just need to estimate the 
${\cal L}(H)$ operator norm
\begin{equation}
\left \|A^{\frac{\kappa}{2}}\ue^{-\beta {\aot}}\right \|_{{\cal L}(H)}\
\leq \sup_{0\leq x<\infty} x^\kappa \ue^{-\beta x} =\left(
\frac{\kappa}{\beta} \right)^\kappa \ue^{-\kappa} = l_\kappa(\beta).
\end{equation}
This concludes the proof of Proposition~2.1.\\

Next, we state and present the proof of the main theorem. 
The steps of the proof are
made in a formal way, however, they can be justified rigorously by 
establishing them first for a Galerkin approximation system and using
the  usual Aubin compactness theorem to pass to the
limit (see, e.g. \cite{lions,constantinfoias,temam}). 
Furthermore, we do not assume that the initial condition $u_0$
is
entire; it is  only assumed to be square integrable, although it will
become
entire for any $t>0$. This is why
in estimating $L^2$ norms of the solution with exponential weights we
have to stay clear of $t=0$.\\

\textbf{Theorem 2.1}\,\, {\it Let $u_0 \in H$, fix $T>0$ and let
$f(t) = f(.\,,t)$ be an entire function with respect to the spatial
variable $x$. Then for every $n=0,1,2\ldots$ there exist constants $C_n,
\bar C_n, K_n$ and $\bar K_n$ which depend on $|u_0|$, $\mu ,
T,\sigma$ and on the norm
\begin{equation}
\int^T_0\left |\ue^{(n-1)\sigma {\aot}}f(s)\right |^2 ds ,
\end{equation}
moreover there exists integers $p_n, q_n\geq 1$ such that
\begin{equation}
\left |\ue^{n \sigma {\aot}}u(t)\right |^2 \leq
\frac{K_n}{t^{p_n}}+C_n, \quad \text{for all}\;\ t \in (0,T]\
\label{eight}\end{equation}
and
\begin{equation}
\int^T_t|\ue^{(n+1) \sigma {\aot}}u(s)|^2 ds\leq \frac{\bar
  K_n}{t^{q_n}}+\bar C_n, \quad \text{for all}  \;t \in (0,T],
\label{nine}
\end{equation}
where $u(t)$ is the solution of \eqref{expns}-\eqref{init}.
}

\textbf{Corollary 2.1} \,\, {\it Let $u_0 \in H, T>0$ and let $f(x,t)$ be an
entire function with respect to the spatial variable $x$ such that for
every $M \geq 0$ we have $\int^T_0 \left |\ue^{M \sigma {\aot}}f(s)\right |^2 ds 
<\infty$. Then, the solution $u(t)$ of \eqref{expns}-\eqref{init} is an entire function with respect 
to the spatial variable for all  $t \in (0,T]$, and satisfies the 
estimates \eqref{eight} and \eqref{nine} in Theorem 2.1 for any $n=1,2,\ldots$.
}\\

\textit{Proof of  Corollary 2.1} \,\, Consider the Fourier series
representation 
\begin{equation}
u(x,t) = \sum_k \ue ^{\ui k\cdot x} \, \hat u(k,t).
\label{fourieru}
\end{equation}
From  \eqref{eight} and Parseval's theorem, we have, for any $n =
0,1,2 \ldots$
\begin{equation}
\sum_k \ue ^{2n \sigma|k|} \left | \hat u(k,t) \right|^2 < \infty,
\label{bidule}
\end{equation}
for $t>0$.
In \eqref{fourieru} we change $x$ to a complex location $z = x+ \ui y$
and obtain
\begin{eqnarray}
u(x+\ui y,t) &=& \sum_k \ue ^{\ui k\cdot (x +\ui y)} \, \hat u(k,t) \\
       &=& \sum_k \left[\ue ^{\ui k\cdot x} \, \ue ^{-|k|}\right]\, 
\left[\ue ^ { -k\cdot y +|k|}\, \hat u(k,t)\right].
\label{plusminus}
\end{eqnarray}
For any $n = 1, 2, \ldots$, the series \eqref{plusminus} of complex
analytic functions  converges uniformly in the strip $|y|+1 \le
\sigma n$. This is because the sum in \eqref{plusminus} is shown to be
bounded, for any $y$,  by use
of the Cauchy--Schwarz inequality applied to the two bracketed 
expressions and use of \eqref{bidule} with $|y|+1 \le \sigma n$. Hence the
Fourier series representation converges in the whole complex domain.
This concludes the  proof of the entire character of the solution with
respect to the spatial variables.\\
\textit{Remark} \,\, This corollary just expresses the most
obvious part of the Paley--Wiener Theorem.\\

\textit{Proof of Theorem 2.1} \,\, The proof of the theorem
proceeds by mathematical induction.
 
\textit{Step $n=0$}\,\, We prove the statement of the
theorem for $n=0$. We take the inner product of \eqref{expns} with $u$ and use
the fact that $(B(u,u),u)=0$ to obtain (when there is no ambiguity
we shall henceforth frequently denote $u(t)$ by $u$)
\begin{eqnarray}
\frac{1}{2}\frac{d}{dt}|u|^2+ \mu \;\left |\ue^{\sigma
{\aot}}u\right |^2&=&(f,u)=(\ue^{-\sigma {\aot}}f ,\ue^{\sigma
{\aot}}u) \\ \nonumber
&\leq& \left |\ue^{-\sigma {\aot}}f \right |\;\left |\ue^{\sigma
  {\aot}}u\right | 
\nonumber \\
&\leq& \frac{|\ue^{-\sigma {\aot}}f |^2}{2\mu}+\frac{\mu }{2}|\ue^{\sigma {\aot}}u|^2 ,
\end{eqnarray}
where Young's inequality has been used to obtain the third line.
Therefore
\begin{equation}
\frac{d}{dt}|u|^2+ \mu \;\left |\ue^{\sigma {\aot}}u\right |^2 \leq \frac{|\ue^{-\sigma {\aot}}f |^2}{\mu}.
\end{equation}
Integrating the above from $0$ to $T$, we obtain
\begin{equation}
|u(t)|^2 + \mu \int^T_0\left |\ue^{\sigma {\aot}}u(s)\right |^2 ds\leq
 C_0 \coloneqq |u_0|^2+ \frac{1}{\mu} \int^T_0\left |\ue^{-\sigma
      {\aot}}f(s)\right|^2 ds.
\label{ten}\end{equation}
Hence
\begin{equation}
|u(t)|^2 \leq C_0,
\label{eleven}
\end{equation}
and
\begin{equation}
\mu \int^T_0\left |\ue^{\sigma {\aot}}u(s)\right |^2 ds \leq C_0.
\label{twelve}
\end{equation}
From \eqref{eleven} and \eqref{twelve} we obtain \eqref{eight} and
\eqref{nine} for the case $n=0$.  Here 
$C_0$ is given by \eqref{ten}, $K_0=0 , \;\bar K_0=0$ and $\bar C_0=C_0$. 
Notice that since $K_0=\bar K_0=0$ there is no need to determine 
the integers $p_0$ and $q_0$ ; however, for the sake of initializing 
the induction process we chose $p_0=q_0=1$.

\textit{Step $n \rightarrow n+1$} 
Assume that \eqref{eight} and \eqref{nine} are true up to $n=m$ and 
we would like to prove them for $n = m + 1$. Let us take the inner product of
\eqref{expns} with $\ue^{2(m + 1)\sigma {\aot}}u $ and obtain
\begin{eqnarray}
\frac{1}{2}\frac{d}{dt}\left |\ue^{(m+1)\sigma {\aot}
}u\right |^2 &+& 
\mu \;\left |\ue^{(m+2)\sigma {\aot}}u\right |^2
\nonumber \\ 
&\leq& \left |(f,\ue^{2(m+1)\sigma {\aot}}u)\right | 
 + \left |(B(u,u),\ue^{2(m+1)\sigma {\aot}}u)\right |
\nonumber \\
&\leq& \left |\ue^{m \sigma {\aot}}f \right |\; \left |\ue^{(m+2)\sigma
  {\aot}}u\right | 
+ \left |\ue^{m \sigma {\aot}}B(u,u)\right
|\;\left |\ue^{(m+2)\sigma {\aot}}u\right |.
\nonumber 
\end{eqnarray}

Now we use Proposition~2.1 to majorize the previous expression by 
\begin{equation}
\leq \left |  \ue ^{ m \sigma \aot } f \right |  
\left |  \ue ^{(m+2) \sigma \aot }  u  \right | + 
C_{\rm A}  \left( l_{\kappa } (\beta ) \right)^\sfa
\left |  \ue ^{ m \sigma  \aot }  u \right |^{2-\sfa}    
\left |  \ue ^{ (m+2) \sigma \aot }  
u  \right |^{1+\sfa}  .
\end{equation}
By Young's inequality we have
\begin{equation}
\begin{split}
&
C_{\rm A} [l_{\kappa } (\beta ) ]^\sfa
\left |  \ue ^{ m \sigma  \aot }  u \right |^{2-\sfa}    
\left |  \ue ^{ (m+2) \sigma \aot }  
u  \right |^{1+\sfa}  \leq
\\
&
\frac{2\kappa - 5 }{4 \kappa }  \, 
\mu ^{-\sfd }  \,  C_{\rm A}^{2\sff }
\, \left( l_{\kappa } (\beta ) \right)^{2\sfa\sff} \, 
\left( \frac{2\kappa + 5}{\kappa } \right)^\sfd \,
\left |  \ue ^{ m \sigma  \aot }  u \right |^{2 \sfe }    + 
\frac{\mu }{4} \, \left |  \ue ^{ (m+2) \sigma \aot }  u  \right |^2 .
\end{split}
\end{equation}
It follows that
\begin{equation}
\begin{split}
&
\frac{d}{dt}  \left | \ue ^{(m+1) \sigma \aot }  u  \right |^2 + 
\mu  \left | \ue ^{(m+2) \sigma \aot }  u  \right |^2 \leq 
\frac{2}{\mu }  \left | \ue ^{ m \sigma \aot } f \right |^2 
\\
&
+\frac{2\kappa - 5}{2 \kappa }  \, 
\mu ^{-\sfd}  \,  C_{\rm A}^{2\sff }
\, \left( l_{\kappa } (\beta ) \right)^{2\sfa\sff} \, 
\left( \frac{2\kappa + 5}{\kappa } \right)^\sfd \,
\left |  \ue ^{ m \sigma  \aot }  u \right |^{2 \sfe}    .
\end{split}
\end{equation}
Now we integrate this inequality on the interval $(s,t) \subset (0,T) $,
obtaining
\begin{equation}
\begin{split}
&
\left | \ue ^{(m+1) \sigma \aot }  u (t)  \right |^2 + 
\mu  \int_s^t 
\left | \ue ^{(m+2) \sigma \aot }  u (s^{\prime } )  \right |^2 \, 
ds^{\prime } 
\\
&
\leq \left | \ue ^{(m+1) \sigma \aot }  u (s)  \right |^2 + 
\frac{2}{\mu }  \int_s^t \left |  \ue ^{ m \sigma \aot } f (s^{\prime } ) 
\right |^2 \, ds^{\prime } +
C' \, \int_s^t
\left |  \ue ^{ m \sigma  \aot }  u (s^{\prime } ) 
\right |^{2 \sfe }   \, ds^{\prime }  
\\
&
\leq \left | \ue ^{(m+1) \sigma \aot }  u (s)  \right |^2 + 
\frac{2}{\mu }  \int_0^T \left |  \ue ^{ m \sigma \aot } f (s^{\prime } ) 
\right |^2 \, ds^{\prime } +
C' \, \int_s^t
\left |  \ue ^{ m \sigma  \aot }  u (s^{\prime } ) 
\right |^{2 \sfe }   \, ds^{\prime }  ,
\end{split}
\label{thirteen}
\end{equation}
where we have set for brevity
\begin{equation}
C' = C (\mu,\beta , \kappa) \coloneqq
\frac{2\kappa - 5 }{2 \kappa }  \, 
\mu ^{-\sfd }  \,  C_{\rm A}
^{2\sff }
\, \left( l_{\kappa } (\beta ) \right)^{2\sfa\sff} \, 
\left( \frac{2\kappa + 5}{\kappa } \right)^\sfd ,
\end{equation}
and where $l_{\kappa } (\beta )$ is given by \eqref{deflkappa}.

Now we come to the point where we use the actual induction assumptions. 
We use \eqref{eight} and the midpoint convexity to estimate the integrand in the last integral: 
\begin{equation}
\left |  \ue ^{ m \sigma  \aot }  u (t) 
\right |^{2 \sfe }  \leq
\left ( \frac{K_m }{t^{p_m } } + C_m \right )
^\sfe  \leq
2^\sff \left ( \frac{K_m }{t^{p_m } }
\right )^\sfe + 
2^\sff C_m^\sfe.
\end{equation}
Whence it follows that
\begin{equation}
\begin{split}
&
C' \int_s^t
\left |  \ue ^{ m \sigma  \aot }  u (s^{\prime } ) 
\right |^{2 \sfe }   \, ds^{\prime } \leq 
2^\sff  C' \, 
\frac{1}{ \sfe p_m - 1} \, 
K_m^\sfe 
\left ( \frac{1}{s^{\sfe p_m  - 1} } - 
\frac{1}{t^{ \sfe p_m  - 1 }}\right
)  
\\
&
+2^\sff  C' \, C_m^\sfe (t-s) \leq 
2^\sff  C' \, 
\frac{1}{\sfe p_m - 1}
K_m^\sfe 
\frac{1}{s^{\sfe p_m   - 1 } }
+ 
2^\sff  C' \, C_m^\sfe (t-s).
\end{split}
\label{missing}
\end{equation}
Discarding the positive term $\mu  \int_s^t 
\left | \ue ^{(m+2) \sigma \aot }  u (s^{\prime } )  \right |^2  \, 
ds^{\prime }  $ in \eqref{thirteen}, we obtain from \eqref{thirteen}
and \eqref{missing}
\begin{equation}
\begin{split}
&
\left | \ue ^{(m+1) \sigma \aot }  u (t)  \right |^2 \leq
\left | \ue ^{(m+1) \sigma \aot }  u (s)  \right |^2 + 
\frac{2}{\mu }  \int_0^T \left |  \ue ^{ m \sigma \aot } f (s^{\prime } ) 
\right |^2 \, ds^{\prime } 
\\
&
+2^\sff  C' \ 
\frac{1}{\sfe p_m  - 1}
K_m^\sfe  
\frac{1}{s^{ \sfe p_m  - 1} } 
+ 
2^\sff  C' \, C_m^\sfe (t-s) .
\end{split}
\end{equation}
Integrating this inequality with respect to $s$ over $(t/2 ,\, t)$ we get
\begin{equation}
\begin{split}
&
\left | \ue ^{(m+1) \sigma \aot }  u (t)  \right |^2 \leq \frac{2}{t} \, 
\int_{t/2}^t  \left | \ue ^{(m+1) \sigma \aot }  u (s)  \right |^2 \, ds + 
\frac{2}{\mu }  \int_0^T \left |  \ue ^{ m \sigma \aot } f (s^{\prime } ) 
\right |^2 \, ds^{\prime } 
\\
&
+\frac{2^\sff  \, C' \, K_m^\sfe }{
(\sfe p_m  - 1)( \sfe p_m 
 - 2 )}
 \left( \frac{2}{t} \right)^{\sfe p_m  - 1} 
+ 
2^\sff  C' \, C_m^\sfe \frac{t}{4}  .
\end{split}
\end{equation}
Note that $p_m\ge1$ implies that 
\begin{equation}
\sfe p_m  - 2  > 0 .
 \end{equation}
By using \eqref{nine}, we have
\begin{equation}
\begin{split}
&
\left | \ue ^{(m+1) \sigma \aot }  u (t)  \right |^2 \leq 
\overline{K} _m \left( \frac{2}{t} \right)^{q_m + 1} + \overline{C} _m \frac{2}{t} +  
\frac{2}{\mu }  \int_0^T \left |  \ue ^{ m \sigma \aot } f (s^{\prime } ) 
\right |^2 \, ds^{\prime } 
\\
&
+\frac{2^\sff  \, C' \, K_m^\sfe }{
( \sfe p_m - 1)( 
\sfe p_m - 2 )}
 \left( \frac{2}{t} \right)^{ \sfe p_m  - 1} 
+ 
2^\sff  C' \, C_m^\sfe \frac{T}{4}  .
\end{split}
\end{equation}
From this relation follows that \eqref{eight} holds for $m+1$ with
\begin{eqnarray}
&&p_{m+1} = \max \left\{  \sfe p_m - 1 ,\, q_m +
  1 \right \},
\label{pmp1}\\
&&q_{m+1} = \max \left \{\sfe  p_m - 1 ,\, q_m + 1
\right \}.
\label{qmp1}
\end{eqnarray}

By the induction assumption we use \eqref{thirteen} to estimate
\begin{equation}
\int^T_{t/2} \left |\ue^{(m+1)\sigma {\aot}}u(s)\right |^2 ds \leq \frac{2^{q_m}\bar K_m}{t^{q_m}}+\bar C_m .
\end{equation}
From this estimate and the above we conclude the existence of the constants $K_{m+1} ,C_{m+1}$ and the integer $p_{m+1}$ such that \eqref{eight} holds for 
$n = m+1$.
Using the estimate that we have just established in \eqref{eight} for $n = m+1$,
and substituting this in \eqref{thirteen}, we immediately obtain the
estimate \eqref{nine} for $n = m+1$. This concludes the proof of Theorem~2.1.

\section{Rate of decay of the Fourier coefficients}
\label{s:decay}

The purpose of this section is to specify the behavior of various
constants appearing in the preceding section to obtain the rate of
decay with the wavenumber of the Fourier coefficients $\hat u(k,t)$ for
$t>0$. We again consider the 3D case in the periodic domain.  Since
the decay may depend on the rate of decay of the Fourier transform of
the forcing term $f(x,t)$, for simplicity we assume zero external
forcing, which we expect to behave as the case with sufficiently rapidly
decaying forcing. The adaptation to sufficiently regular forced cases, 
for example a trigonometric polynomial, is similar but more technical.\footnote{It is
conceivable 
that the results can be extended to  forces entire in the space
variables whose Fourier transforms decrease faster than $\ue ^{-C|k|
  \ln|k|}$ with sufficiently large $C$.}  Furthermore, it is enough to
prove the decay result up to a  time $T$ such that 
 $1/N \coloneqq TU/L <1$, where $L$ and $U$ are a typical length scale and velocity of
the initial data.  Extending the results to later times is easy (by
propagation of regularity).

We shall show that the bound for the square of the 
 $L^2$ norm of the velocity weighted by $\ue ^{n\sigma \aot}$ is a
 double exponential in $n$. Specifically, we have\\

\textbf{Theorem 3.1} \,\,{\it Let $u(t)$ be the solution of
\eqref{expns}-\eqref{init} in $[0,\,T]$ with $f=0$ and  $ 0 < T <
L/U$. Then
for every $\kappa>3$ and $\delta > 0 $,
there exists a   number $\Lambda$, depending on
$\delta$ and $\kappa$,
such that, for all integer $n \ge 0$
\begin{eqnarray}
\left |  \ue ^{ n\sigma \aot}  u (t)  \right |^2 &\leq&
\left( \frac{\Lambda L}{U t} \right)^{a_n}   
, \qquad t \in (0,T] ,
\label{e:estimate9}\\
\int_t^T \left | \ue ^{ (n+1) \sigma \aot }  u (s)  \right|^2 \, ds &\leq&
\left( \frac{\Lambda L}{U t}  \right)^{ a_n} 
, \qquad t \in (0,T] 
\label{e:estimate10}\\
 {\it where} \,\,\,\,\,\,a_n &=& \left(( 1 + \delta ) \, \frac{4\kappa - 5}{2\kappa - 5}\right)^n.
\label{defan}
\end{eqnarray}
}

\textbf{Corollary 3.1} \,\, {\it For any $t > 0$ the function $u(t)$
of \eqref{expns}-\eqref{init} is an entire function in the space variable and its (spatial) Fourier coefficients  
tend to zero in the following faster-than-exponential way: there exists
a constant $\Lambda$ such that, for any $0<\varepsilon<1$, we have
\begin{equation}
\vert \hat{u} (k,t) \vert  \leq  \ue ^{- \frac{\sigma (1 - \varepsilon ) }{
\beta _{\kappa , \delta } } 
\vert k \vert \ln \vert k \vert } , \quad \hbox{ for all}\quad \vert k \vert \geq \left(
\sqrt{\frac{\Lambda L}{U t} } \right)^{
\frac{\beta _{\kappa , \delta } }{\varepsilon \sigma } } , 
\end{equation}
where 
\begin{equation}
\beta _{\kappa , \delta } = 
\ln \left( (1 + \delta ) \, \frac{4\kappa - 5}{2\kappa - 5} \right) . 
\end{equation}
}

\textit{Proof of Corollary 3.1} \,\,Since we are dealing with a
Fourier series, the modulus of any Fourier coefficient of the function
 $\ue ^{ (n+1) \sigma \aot }  u(t)$ cannot exceed its
 $L^2$ norm, hence it is bounded by \eqref{e:estimate9}. Thus, discarding a
factor
$(2\pi)^{-3/2} <1$, we have for all
$k $ and $n$ 
\begin{equation}
\vert \hat{u} (k,t) \vert  \leq    
\ue ^{-n\sigma \vert k \vert } 
\left( \sqrt{\frac{\Lambda L}{U t } }\right)^{\left( (1 + \delta ) \, 
\sfe \right)^n }  = \exp \Bigl( \ln \sqrt{\frac{\Lambda L}{U t} }  
\ue ^{n \beta _{\kappa , \delta } }  - n \sigma \vert k \vert \Bigr) ,
\end{equation}
where $\sfe$ is defined in \eqref{defsfag}.
Now choosing 
\begin{equation}
\ln \vert k \vert \geq \frac{1}{\varepsilon } \frac{
\beta _{\kappa , \delta } }{\sigma } \ln 
\sqrt{\frac{\Lambda L}{U t} } , 
\end{equation}
we obtain with $ n = \ln \vert k \vert / \beta _{\kappa , \delta } $ the following estimate
\begin{equation}
\vert \hat{u} (k,t) \vert  \leq 
 \exp \biggl[   - \frac{\sigma }{\beta _{\kappa , \delta } }  
\vert k \vert \ln \vert k \vert 
 \Bigl( 1 -  \frac{\beta _{\kappa , \delta }  \ln 
\sqrt{\Lambda L/U t} }{\sigma  \ln \vert k \vert }   \Bigr)
 \biggr] \leq \ue ^{- \frac{\sigma (1 - \varepsilon ) }{\beta _{\kappa , \delta } } 
\vert k \vert \ln \vert k \vert } .
\end{equation}

\textit{Remark 3.1} \,\, Since $\varepsilon$ and $\delta $ can be chosen
arbitrarily small and $\kappa$ arbitrarily large, Corollary~3.1
implies that, in terms of the dimensionless wavenumber
$\kt = 2\sigma k$, the Fourier amplitude has a  bound (at high
enough $\kt$) of the form $\ue ^{- C\kt \ln \kt}$ for any $C <1/(2
\ln 2)$. We shall see  that the upper bound for the constant $C$ can
probably be improved to $1/\ln 2$.

\textit{Proof of Theorem 3.1} \,\,
We proceed again by induction. We assume that the following inequalities hold
\begin{equation}
\label{e:ass1}
\left | \ue^{n \sigma \aot }  u  (t) \right |^2
\leq     
\frac{K_n }{t^{a_n } }  ,
\end{equation}
and 
\begin{equation}
\label{e:ass2}
\int_t^T \left | \ue^{(n+1) \sigma \aot }  u  (s) \right | ^2\, ds
\leq  
\frac{K_n }{t^{a_n} }  ,
\end{equation}
where $K_n$ and $a_n \ge 1$ are still to be determined.
Starting from  expNSE \eqref{expns}, we take the inner product with
$\ue^{2 (n+1) \sigma \aot} u $. Then we obtain from Proposition~2.1
with $\alpha = n\sigma $ and $\beta = 2 \sigma $
\begin{equation}
\begin{split}
&
\frac{1}{2} \frac{d}{dt}  \left | \ue^{(n+1) \sigma  \aot}  u  \right |^2 + 
\mu  \left | \ue^{(n+2) \sigma \aot }  u  \right |^2 \leq 
 \left | \bigl( B ( u , u ) , \ue^{2 (n+1) \sigma \aot }  u \Bigr)
 \right | 
\\
&\leq 
\left |  \ue^{ n \sigma \aot }  B ( u , u ) \right |\,
\left |  \ue^{ (n+2) \sigma \aot }  u  \right |  \leq
C_{\rm A}  \bigl( l_{\kappa } (\beta ) \bigr)^\sfa
\left |  \ue^{ n \sigma  \aot }  u \right | ^{2-\sfa}    
\left |  \ue^{ (n+2) \sigma \aot }  
u  \right |^{1+\sfa}   
\\
&
\leq\frac{2\kappa - 5 }{4 \kappa }  \, 
\mu ^{-\sfd }  \,  C_{\rm A }^{2\sff}
\, \Bigl( l_{\kappa } (\beta ) \Bigr)^{2\sfa\sff} \, 
\left( {1+\sfa} \right)^\sfd \,
\left |  \ue^{ n \sigma  \aot }  u \right |^{2 \sfe }    + 
\frac{\mu }{2} \, \left |  \ue^{ (n+2) \sigma \aot }  u  \right |^2 .
\end{split}
\end{equation}
Then it follows that
\begin{equation}
\begin{split}
&
\frac{d}{dt}  \left | \ue^{(n+1) \sigma \aot }  u  \right |^2 + 
\mu  \left | \ue^{(n+2) \sigma \aot }  u  \right | ^2
\\
&
\leq 
\frac{2\kappa - 5 }{2 \kappa }  \, 
\mu ^{-\sfd }  \,  C_{\rm A }^{2\sff }
\, \Bigl( l_{\kappa } (\beta ) \Bigr)^{2\sfa\sff} \, 
\left( \frac{2\kappa + 5}{2 \kappa } \right)^\sfd \,
\left |  \ue^{ n \sigma  \aot }  u \right |^{2 \sfe } .   
\end{split}
\end{equation}
By using the induction assumption we obtain
\begin{equation}
\label{e:estimate1}
\frac{d}{dt}  \left | \ue^{(n+1) \sigma \aot }  u  \right |^2 + 
\mu  \left | \ue^{(n+2) \sigma \aot }  u  \right |^2 
\leq 
C^{\prime\prime}  \,
\left( \frac{K_n }{t^{a_n } }  \right)^\sfe,    
\end{equation}
where we have set 
\begin{equation}
C^{\prime\prime} = \frac{2\kappa - 5 }{2 \kappa }  \, 
\mu ^{-\sfd}  \,  C_{\rm A }^{2\sff }
\, \Bigl( l_{\kappa } (\beta ) \Bigr)^{2\sfa\sff} \, 
\left( \frac{2\kappa + 5}{2 \kappa } \right)^\sfd.
\end{equation}
Renaming the time variable in \eqref{e:estimate1} from $t$ to $s'$ and
integrating over $s'$ from $s$ to $t$ \hbox{(with $0 <s<t \le T$)} we obtain
\begin{equation}
\begin{split}
&
\left | \ue^{(n+1) \sigma \aot }  u  (t) \right |^2 + 
\mu  \int_s^t \left | \ue^{(n+2) \sigma \aot }  u (s^{\prime } )  \right | ^2
d s^{\prime }
\\
&
\leq 
\frac{1}{a_n \sfe  - 1}
C^{\prime\prime}  \, K_n^\sfe  
\left(  \frac{1}{s^{a_n\sfe  - 1} }  -  \frac{1}{t^{a_n \sfe - 1} } \right) + 
\left | \ue^{(n+1) \sigma \aot }  u  (s) \right | ^2 
\\
&
\leq \frac{1}{a_n \sfe  - 1}
C^{\prime\prime}  \, K_n^\sfe   
\frac{1}{s^{a_n \sfe  - 1} }   + 
\left | \ue^{(n+1) \sigma \aot }  u  (s) \right |^2.
\end{split}
\end{equation}
Omitting the positive integral term on the left-hand side of the inequality we obtain
\begin{equation}
\left | \ue^{(n+1) \sigma \aot }  u  (t) \right |^2 
\leq 
\frac{1}{a_n \sfe  - 1}
C^{\prime\prime}  \, K_n^\sfe  
\frac{1}{s^{a_n \sfe  - 1} }   + 
\left | \ue^{(n+1) \sigma \aot }  u  (s) \right |^2.
\end{equation}
Choosing $ 1 < \gamma \leq N^\delta = \left(L/UT\right)^{\delta } $ and 
integrating over $s$ from $t/\gamma$ to $t$ we obtain
\begin{equation}
\label{e:onehalf}
\begin{split}
&
(\gamma - 1) \frac{t}{\gamma }
\left | \ue^{(n+1) \sigma \aot }  u  (t) \right |^2 
\leq 
\frac{1}{a_n \sfe  - 1}
\frac{1}{a_n \sfe  - 2}
C^{\prime\prime}  \, K_n^\sfe   
\biggl\{ \Bigl( \frac{\gamma}{t} \Bigr) ^{a_n\sfe  - 2}   - 
\Bigl( \frac{1}{t} \Bigr) ^{a_n \sfe  - 2}   \biggr\}
\\
&
+ 
\int_{t/\gamma}^t 
\left | \ue^{(n+1) \sigma \aot }  u  (s) \right |^2 \, ds \leq
\frac{1}{a_n \sfe  - 1}
\frac{1}{a_n \sfe  - 2}
C^{\prime\prime}  \, K_n^\sfe   
\Bigl( \frac{\gamma}{t} \Bigr) ^{a_n \sfe  - 2}   +
K_n \Bigl( \frac{\gamma}{t} \Bigr)^{a_n }  ,
\end{split}
\end{equation}
where we have used the induction assumption \eqref{e:ass2}. We obtain thus the following 
estimate
\begin{equation}
\label{e:estimate2}
\left | \ue^{(n+1) \sigma \aot }  u  (t) \right | ^2
\leq 
\frac{1}{\gamma - 1 } \, \frac{1}{a_n \sfe  - 1}
\frac{1}{a_n \sfe  - 2}
C^{\prime\prime}  \, K_n^\sfe  
\Bigl( \frac{\gamma}{t} \Bigr) ^{a_n \sfe  - 1}   +
\frac{1}{\gamma - 1} \, K_n \Bigl( \frac{\gamma}{t} \Bigr)^{a_n + 1} ,
\end{equation}
which holds for every 
$ 0 < t \leq T$.

To estimate $\int_t^T\left| \ue ^ {(n+2)\sigma \aot} u(s)\right|^2\,ds$ we integrate \eqref{e:estimate1} from $t$ to $T$:
\begin{equation}
\begin{split}
&
\left | \ue^{(n+1) \sigma \aot }  u (T)  \right |^2 + 
\mu  \int_t^T \left | \ue^{(n+2) \sigma \aot }  u (s)  \right |^2 
\, ds
\\
&
 \leq C^{\prime\prime}  \,
\frac{1}{a_n\sfe  - 1}
C^{\prime\prime}  \, K_n^\sfe   
\frac{1}{t^{a_n\sfe  - 1} }  
+ \left | \ue^{(n+1) \sigma \aot }  u (t)  \right | ^2.
\end{split}
\end{equation}
Omitting the first term on the right-hand side  and using \eqref{e:estimate2} we obtain
\begin{equation}
\label{e:estimate3}
\begin{split}
&
\mu  \int_t^T \left | \ue^{(n+2) \sigma \aot }  u (s)  \right | ^2
\, ds \leq 
C^{\prime\prime}  \,
\frac{1}{a_n \sfe  - 1}
C^{\prime\prime}  \, K_n^\sfe   
\frac{1}{t^{a_n \sfe  - 1} }  
\\
&
+\frac{1}{\gamma - 1} \, \frac{1}{a_n \sfe  - 1}
\frac{1}{a_n \sfe  - 2}
C^{\prime\prime}  \, K_n^\sfe   
\Bigl( \frac{\gamma}{t} \Bigr) ^{a_n \sfe  - 1}   +
\frac{1}{\gamma - 1} \, 
K_n \Bigl( \frac{\gamma}{t} \Bigr)^{a_n + 1} .
\end{split}
\end{equation}
We conclude that since $a_n \geq 1 $ and 
\begin{equation}
a_n + 1 \leq a_n \Bigl(2 + \frac{5}{2\kappa - 5} \Bigr) 
= a_n \sfe,
\end{equation}
for a suitable constant $E>0$ we have
\begin{equation}
\label{e:estimate3b}
\left | \ue^{(n+1) \sigma \aot }  u  (t) \right | ^2
\leq  
E K_n ^\sfe  
\left( \frac{\gamma }{t} \right)^{a_n\sfe  }  
\end{equation}
and 
\begin{equation}
\label{e:estimate4}
\int_t^T \left | \ue^{(n+2) \sigma \aot }  u  (s) \right | ^2
\leq  
E K_n ^\sfe   
\left( \frac{\gamma }{t} \right)^{a_n \sfe  } . 
\end{equation}
Since,  for $t \leq T $, 
\begin{equation}
\frac{\gamma} {t} \leq \frac{N^\delta}{t} = \left(\frac{L}{U
  T}\right)^\delta \frac{1}{t}\leq \left(\frac{L}{Ut}\right)^{1+\delta} \frac{U}{L}, \nonumber
\end{equation} 
it follows that
\begin{equation}
\label{e:estimate3c}
\left | \ue^{(n+1) \sigma \aot }  u  (t) \right | ^2
\leq  
E K_n ^\sfe  
\left( \frac{L}{U t} \right)^{a_n (1+\delta ) \sfe }
\left(\frac{U}{L}\right)^{a_n \sfe}
\end{equation}
and 
\begin{equation}
\label{e:estimate4b}
\int_t^T \left | \ue^{(n+2) \sigma \aot }  u  (s) \right | ^2
\leq  
E K_n ^\sfe   
\left( \frac{L}{U t} \right)^{a_n (1+\delta ) \sfe  }\left(\frac{U}{L}\right)^{a_n \sfe} . 
\end{equation}
This finishes the induction step. 

From the above follows that we can take
\begin{equation}
a_{n+1} = a_n (1+\delta ) \frac{4\kappa - 5}{2\kappa - 5} , \qquad
K_{n+1} = E K_n^\sfe .
\end{equation}
Note that in the induction step we use the assumption that $a_n \geq 1$. This fixes the value of 
$a_0 = 1$. The solution of the recursion relations is given by
\begin{equation}
a_n = \left( (1 + \delta ) \, \frac{4\kappa - 5}{2\kappa - 5} \right)^n ,
\qquad
\ln K_{n} = \ln E \, \frac{a_n - 1}{ (1 + \delta ) \, 
\frac{4\kappa - 5}{2\kappa - 5} - 1 } + a_n\, \ln K_0 .
\end{equation}
Finally, choosing a sufficiently large number $\Lambda$ we get the desired
estimates \eqref{e:estimate9} and \eqref{e:estimate10}. This concludes
the proof of Theorem~3.1.

\section{Remarks and extensions for the main results}
\label{s:extensions}

Although our main theorems are stated for the case of the 3D expNSE,
their statements and proofs are easily extended \emph{mutatis mutandis}
to arbitrary space dimensions $d$: with exponential dissipation for any $d$ the solution is entire
in the space variables and the decay of Fourier coefficients is
bounded by $\exp(-C |\kt| \ln|\kt|)$ for any $C< C_\star =1/(2\ln 2)$. Some of the
intermediate steps in the proof, such as the formulation of Agmon's
inequality, change with $d$ but not the result about the constant $1/(2\ln 2)$.

We can also easily change the functional form of the dissipation.\footnote{Note that  the 
proof of Proposition~2.1 and Theorem~2.1 holds \emph{mutatis mutandis} if 
we replace, in the argument of the exponential, $|k|$ by  a subadditive
function of $|k|$ subject to some mild conditions, such as
$|k|^{\alpha}$ with $0<\alpha <1$.} 
One instance is a dissipation operator ${\cal D}$ with a Fourier symbol $\ue
^{2\sigma |k|^{\alpha}}$ with $0<\alpha<1$.  One can prove that the solution in this case satisfies
\begin{equation}
\sum_{k \neq 0} \ue ^{2n|k|^\alpha} \left|\hat u(k,t)\right|^2 \leq
\frac{K_n}{t^{p_n}} +C_n,
\label{alphamod}
\end{equation}
for all $t \in (0,T]$ and for all $n$. Hence the solution in this case
  belongs to $C^\infty$ but is not necessarily an entire function. In
  fact
it belongs to the Gevrey class $G_{1/\alpha}$. Gevrey regularity with
$0<\alpha <1$ does not even imply analyticity.\footnote{The special
  class
when $\alpha =1$ of analytic functions is considered by some authors as
one of the Gevrey classes \cite{foiastemam,ferrarititi}.} Actually,
 with such a
dissipation,
the solutions are analytic even when $\alpha<1$.  We shall return to this case of dual Gevrey regularity and analyticity
in Section~\ref{ss:heuristics}.   

Next, consider the case  $\alpha >1$. The dissipation has a lower
bound of the ordinary exponential type, so that the entire character of the
solution is easily established.  However, for $\alpha>1$ the bound $\exp(-C
|k| \ln|k|)$ can be improved in its functional form, as we shall see in
Section~\ref{ss:heuristics}.  

Obviously, the results of
Sections~\ref{s:entire} and \ref{s:decay} do still hold if we change the
functional form of the Fourier symbol of the dissipation at low wavenumbers
$|k|$ while keeping its exponential growth at high wavenumber. One
particularly interesting instance, to which we shall come back in the
next Section
on the Burgers equation and in the Conclusion, is ``cosh dissipation'', namely
a Fourier symbol $-\mu(1-\cosh (k/\kd))$ with $\mu>0$. The dissipation rate at
wavenumber much smaller than $\kd$ is then $\nu |k|^2$ with $\nu =
\mu/(2\kd^2)$, just as for the ordinary Navier--Stokes equation.

It is worth mentioning that the key results of Sections~\ref{s:entire}
and
\ref{s:decay} still hold when the problem is formulated in the whole 
space $\mathbb{R}^d$ rather than with periodicity conditions.  Similarly they should hold on
the sphere $S^2$, a case for which spherical harmonics can be used
(see \cite{CRT}).

Of course the result on the entire character of the solution, 
when exponential dissipation is assumed, holds for a large class of partial
differential equations. Besides the exponential modification of the  Navier--Stokes equations it
applies  to similar modifications, for example,  of the magnetohydrodynamical equations 
and  of the complex Ginzburg--Landau equation 
\begin{equation}
\frac{\partial u}{\partial t}- \alpha\frac{\partial ^2 u}{\partial x^2} + \beta u +\gamma |u|^2 u=0,
\label{}\end{equation}
where
\begin{equation} 
{\rm Re}\:\ \alpha >0\:\ ,\:\ {\rm Re}\:\ \gamma>0.
\label{}\end{equation}
The main idea would be  in proving the analogue of Proposition~2.1 for
the corresponding nonlinear terms in the underlying equations
following our proof combined with ideas presented in 
\cite{ferrarititi} and  \cite{doelmantiti}.

\section{The case of the 1D Burgers equation}
\label{s:burgers}

The (unforced) one-dimensional Burgers equation with modified
dissipation reads:
\begin{eqnarray}
&&\frac{\partial u}{\partial t} + u\frac{\partial u}{\partial x}= - \mu {\cal D} u,\label{burgers}\\[1.ex]
&&u(x,0)=u_0(x). \label{initburgers}
\end{eqnarray}
We shall mostly consider the case of the cosh Burgers equation when ${\cal D}$
has the Fourier symbol $-\mu(1-\cosh(k/\kd))$. Since the cosh Burgers
equation is much simpler than expNSE we can expect to obtain stronger results
or, at least, good evidence in favor of stronger conjectures.

Let us observe that the cosh Burgers equation can be rewritten
in the complexified space of analytic functions of $z\coloneqq x+ \ui y$ as 
\begin{equation}
\frac{\partial u(z,t)}{\partial t} 
+ u(z,t)\frac{\partial u(z,t)}{\partial x}= 
-  \frac{\mu}{2}\left[u(z+\ui/\kd, t)+u(z-\ui/\kd, t) -2u(z,t)\right].
\label{imdiffburgers}
\end{equation}
This is the ordinary Burgers equation with the dissipative Laplacian
replaced by its centered second-order finite difference approximation,
differences being taken in the pure imaginary direction with a mesh
$1/\kd$.
 
As already stated, Corollary~2.1 on the entire character of the solution and
Corollary~3.1 on the bound of the modulus of the Fourier coefficients by
$\exp(-C |\kt| \ln|\kt|)$ for any $C <1/(2\ln 2)$ hold in the same form as
for the expNSE. Of course, if the finite differences were taken in the real
rather than in the pure imaginary direction, the solution would not be entire.
Actually, \eqref{imdiffburgers} relates the values of the velocity on lines
parallel to the real axis shifted by $\pm 1/\kd$ in the imaginary direction.
It thereby provides a kind of \emph{Jacob's Ladder} allowing us to
climb to complex infinity in the imaginary direction. This can be used
to show, at least heuristically, that  the complexified velocity grows
with the imaginary coordinate $y$ as $\exp\left(C 2^{|y|\kd}\right)$. 

Such a heuristic derivation turns out to be equivalent to another
derivation by dominant balance which can be done on the
Fourier-transformed equation, the latter being not limited to cosh
dissipation. Section~\ref{ss:heuristics} is devoted to Fourier space
heuristics for different forms of the dissipation.  For exponential
and cosh dissipation this suggests a leading-order behavior of the
Fourier coefficients for large wavenumber of the form $\exp(-C_\star
|k| \ln|k|)$ with $C_\star =1/\ln 2$, a substantial improvement over
the rigorous bound. Various numerical and semi-numerical results,
discussed in Section~\ref{ss:simulation1} and \ref{ss:simulation2},
support this improved result.

\subsection{Heuristics: a dominant balance approach}
\label{ss:heuristics}

We want to handle dissipation operators ${\cal D}$
with an arbitrary positive Fourier symbol, taken here to be $\ue ^{G(k)}$
where $G(k)$ is a real even function of the wavenumber
$k\in\mathbb{Z}$ which is increasing without bound for $k>0$.  It is then best to rewrite the Burgers equation in terms
of the Fourier coefficients. We set 
\begin{equation}
u(x,t) =\sum_{k\in\mathbb{Z}} \ue ^{\ui kx}\, \hat u(k,t),
\label{deffourier}
\end{equation}
and obtain from \eqref{burgers}
\begin{equation}
\frac{\partial \hat u(k,t)}{\partial t} +\frac{\ui k}{2}\sum_{p+q=k} \hat
u(p,t) \hat u(q,t) = -\mu \ue ^{G(k)} \hat u(k,t).
\label{fourierburgers}
\end{equation}
This is the place where we begin our heuristic analysis of the 
large-wavenumber asymptotics. First, we
drop the time derivative term since it will turn out not to be relevant to
leading order. (A suitable Galilean change of frame may be needed
before this becomes true.) For simplicity we now drop the time variable
completely. The next heuristic step is to balance the moduli of the
two remaining terms, taking 
\begin{equation}
|\hat u(k) | \sim \ue ^{-F(k)},
\label{defFk}
\end{equation}
where $F(k)$ is still to be determined but assumed sufficiently smooth
and the symbol $\sim$ is used here
to connect two functions ``asymptotically equal up to constants and algebraic 
prefactors'' (in other words, asymptotic equality of the logarithms).
The convolution in \eqref{fourierburgers} can be approximated for
large wavenumbers by a continuous wavenumber integral $\sim 
\int \ue ^{-F(p)-F(k-p)}  dp$. Next we evaluate the integral
by steepest descent, assuming that the leading order comes from the critical
point $p= k/2$, where the $p$-derivative of $F(p)+F(k-p)$ obviously
vanishes. This will require that this point be truly a minimum of 
$F(p)+F(k-p)$. Balancing the logarithms of the nonlinear term and
of the dissipative term we obtain the following simple equation for
the function $F(k)$:
\begin{equation}
2F\left(\frac{k}{2}\right) = F(k) -G(k).
\label{FGequation}
\end{equation}
This is a linear first order finite difference equation (in the
variable $\ln k$) which is easily solved for values of the wavenumber
 of the form $k = 2^n$:
\begin{equation}
F(2^n) = 2^n\left[F(1) +\frac{G(2)}{2}+\frac{G(4)}{4}+\dots +\frac{G(2^n)}{2^n}\right]. 
\label{explicitF}
\end{equation}

For exponential dissipation (and for cosh dissipation when $|k|/\kd \gg
1$), we have 
 $G(k) \simeq 2\sigma |k|$ and we obtain
from \eqref{explicitF}, to leading order for large positive $k$
\begin{equation}
F(k) \simeq \frac{1}{\ln 2} \kt \ln \kt; \quad \kt\coloneqq 2\sigma k = \frac{k}{\kd}.
\label{Fforexp}
\end{equation}
If this heuristic result is correct -- and the supporting numerical
evidence is strong as we shall see in  Sections~\ref{ss:simulation1}
 and \ref{ss:simulation2} -- the estimate given by Corollary~3.1 (adapted
to the Burgers case) that $|\hat u(k)| < \ue ^{-C |\kt|\ln |\kt|}$ for
sufficiently large $|\kt|$ and any $C <C_\star =1/(2\ln 2)$  still leaves
room for improvement as to the value of $C_\star$. It can be shown that this
dominant balance argument remains unchanged if we reinsert the
time-derivative term, since its contribution is easily checked to be 
subdominant. Actually, the conjecture that the solution of NSE is entire
with exponential or cosh dissipation was based on precisely this kind of
dominant balance argument, which suggests a faster-than-exponential
decay of the Fourier coefficients.

When $G(k) = 2 \sigma |k|^\alpha$ with $\alpha>1$ we obtain to leading
order
\begin{equation}
F(k) \simeq \frac{2\sigma}{1 -2^{1-\alpha}} \,|k|^\alpha.
\label{Fforalpha}
\end{equation}
This is an even faster decay of the Fourier coefficients than in the
exponential case \eqref{expns}.\footnote{Actually, one can show, for the Burgers
  equation and the NSE that when 
$\alpha>1$ the Fourier coefficients of the solution decay faster
than $\exp (-C|k|^{\alpha -\epsilon})$ for any $\epsilon > 0 $.}

It is easily checked that for $\alpha\ge 1$ the condition of having a
minimum of $F(p)+F(k-p)$ at $p=k/2$ is satisfied. If however we were
to use \eqref{Fforalpha} for $0<\alpha<1$ the condition would not be
satisfied. In this case it is easily shown for the Burgers
 equation and the NSE, by using a variant of the
theory presented in Section~\ref{s:entire},  that the solution is in the Gevrey
class $G_{1/\alpha}$ in the whole space $\mathbb{R}^d$.  It is actually not difficult to show that
the solution is also analytic when $0<\alpha<1$, in a finite strip in
$\mathbb{C}^d$ about the real space $\mathbb{R}^d$. For this it suffices
to adapt to the proof of analyticity given for the  ordinary NSE under 
the condition of some mild regularity. Such regularity is trivially satisfied
with the much stronger dissipation assumed here
\cite{foiastemam}.\footnote{The first results on analyticity, derived in the more complex
  setting of flow with boundaries, were obtained in \cite{masuda}.}
We also found strong numerical evidence for analyticity. 

It is of interest to point out that, although analyticity is a stronger
regularity than Gevrey when $0<\alpha<1$, the Gevrey result implies
a decay of the  form $\exp (-C|k|^\alpha \ln |k|)$,
independently
of the viscosity coefficient $\mu$, whereas analyticity in a finite
strip gives a decay
of the form $\exp(- \eta |k|)$ where $ \eta$ depends on $\mu$ 
\cite{doeringtiti}.

\subsection{Spectral simulation for the Burgers case}
\label{ss:simulation1}

Here we begin our numerical tests on the 1D Burgers equation. So far
we have a significant gap in the value of the constant $C$ appearing
in the $\ue ^{-C |\kt|\ln |\kt|}$ estimate of Fourier coefficient,
between the bounds and a heuristic derivation of the asymptotic behavior.
In this section we shall exclusively consider the case of the unforced
Burgers equation with initial condition $u_0(x) = -\sin x$ and
dissipation with a rate $1-\cosh k$.  (Thus, $\mu =1$ and $\kd =1$.)
The numerical method is however very easily extended to other
functional forms of the dissipation and other initial conditions. The
spectral method is actually quite versatile. Its main drawback will be
discussed at the end of this section.

The standard way of obtaining  a high-orders scheme when numerically
integrating PDE's with (spatial)
periodic boundary conditions is by the (pseudo)-spectral technique
with
the 2/3 rule of alias removal \cite{gottlieborszag}.
The usual  reason this is more precise than finite differences is 
that the truncation errors resulting from the use of a finite number
$N$ of collocation points (and thus a finite number $N/3$ of Fourier
modes) decreases exponentially with $N$ if the
solution is analytic in a strip of width $\delta$ around the real
axis. Indeed
this implies a bound for the Fourier coefficients at high $|k|$ of the form
$|\hat u(k)| < \ue ^{-C|k|}$ for any $C<\delta$. In the present case,
the solution being entire, the bound is even better. 

There are of course sources of error other than spatial Fourier
truncation,
namely rounding errors and temporal discretization errors. Temporal
discretization is a non-trivial problem here because the dissipation
grows exponentially with $|k|$ and thus the characteristic time scale
of high-$|k|$ modes can become exceedingly small. Fortunately, these
modes are basically slaved to the input stemming from nonlinear
interaction of lower-lying modes. It is possible to take advantage
of this to use a slaving technique which bypasses the stiffness of the
equation (a simple instance of this phenomenon is
described in Appendix~B of \cite{FST}). We use here the slaved scheme
\emph{Exponential Time Difference Runge Kutta 4} (ETDRK4) of  
\cite{coxmatthews}
with a time step of $10^{-3}$\;.\footnote{This is far larger than would
have have been permitted without the slaving. Actually it can still
be increased somewhat to $5\times 10^{-3}$ without affecting the
results.}

As to the rounding noise, it is essential to use at least double 
precision since otherwise the faster-than-exponential decrease of the
Fourier coefficients would be swamped by rounding noise beyond
a rather modest wavenumber. Even with double precision, rounding noise
problems start around wavenumber 17, as we shall see. Hence it makes
no sense to use more than, say, 64 collocation points, as we have done.

\begin{figure}
\centerline{%
\includegraphics[scale=0.50,angle=-90]{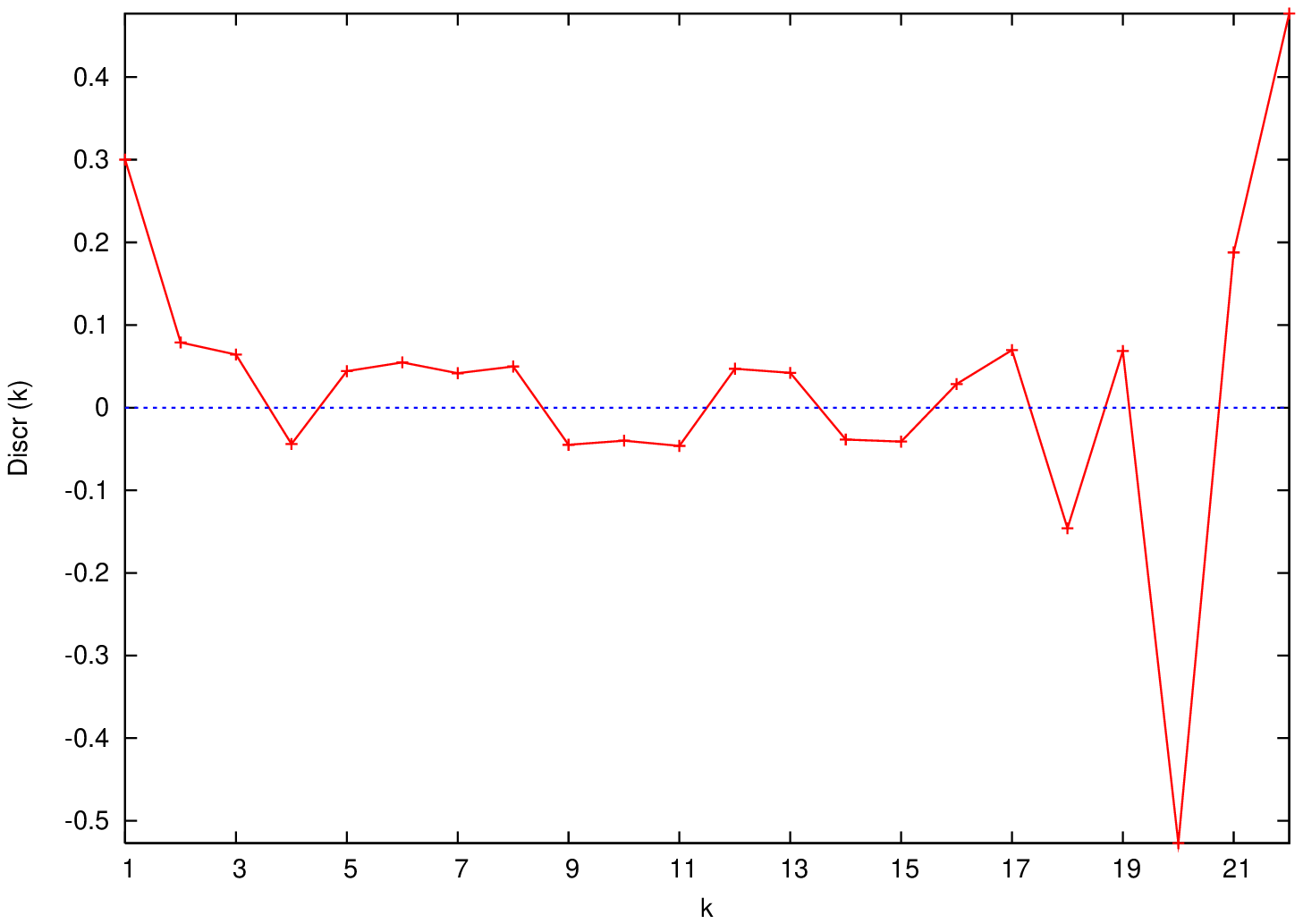}%
}
 \caption{Discrepancy Discr, as predicted by  \eqref{discrepancyspectral}, vs 
wavenumber $k$ for data obtained by spectral simulation in double precision;
rounding noise becomes significant beyond wavenumber 17.}
\label{f:stupid-discrepancy-spectral}
\end{figure}

Fig.~\ref{f:stupid-discrepancy-spectral} shows the discrepancy
\begin{equation}
{\rm Discr}\,(k) \coloneqq -\frac{\ln|\hat u(k,1)|} {|k| \ln |k|} -\frac{1}{\ln 2}, 
\label{discrepancyspectral}
\end{equation}
which, according to heuristic asymptotic theory \eqref{defFk}-\eqref{Fforexp}, 
should converge to $0$ as $|k| \to \infty$. 

It is seen that the discrepancy falls to about 3.5\% of the nominal
value $1/\ln 2$ before getting swamped by
rounding noise around wavenumber 17.

It is actually possible to significantly decrease the discrepancy
by using a better processing of the numerical output, called \textit{asymptotic
extrapolation}, developed recently by van der Hoeven \cite{jorasint} 
and which is related to the theory of transseries
\cite{ecalle,jorisspringer}. The basic idea is to perform on the data a sequence of
transformations which successively strip off the leading and subleading
terms in the asymptotic expansion (here for large $|k|$). Eventually,
the transformed data allow a very simple interpolation (mostly
by a constant). The procedure can be carried out until the transformed
data
become swamped by rounding noise or display lack of asymptoticity, 
whichever occurs first. After the interpolation stage, the successive 
transformations
are undone. This determines the asymptotic expansion of the data
up to a certain order of subdominant terms.  
An elementary introduction to this method
may be found in  \cite{PF06}, from which we shall also borrow the notation for the various
transformations: {\bf I} for ``inverse'', {\bf R} for ``ratio'',
{\bf SR} for ``second ratio'', {\bf D} for ``difference'' and {\bf
  Log} for ``logarithm''.
The choice of the successive transformations is dictated by
various tests which roughly allow to find into which 
broad asymptotic class the data and their transformed versions
fall. 
\begin{figure}
\centerline{%
\includegraphics[scale=0.50,angle=-90]{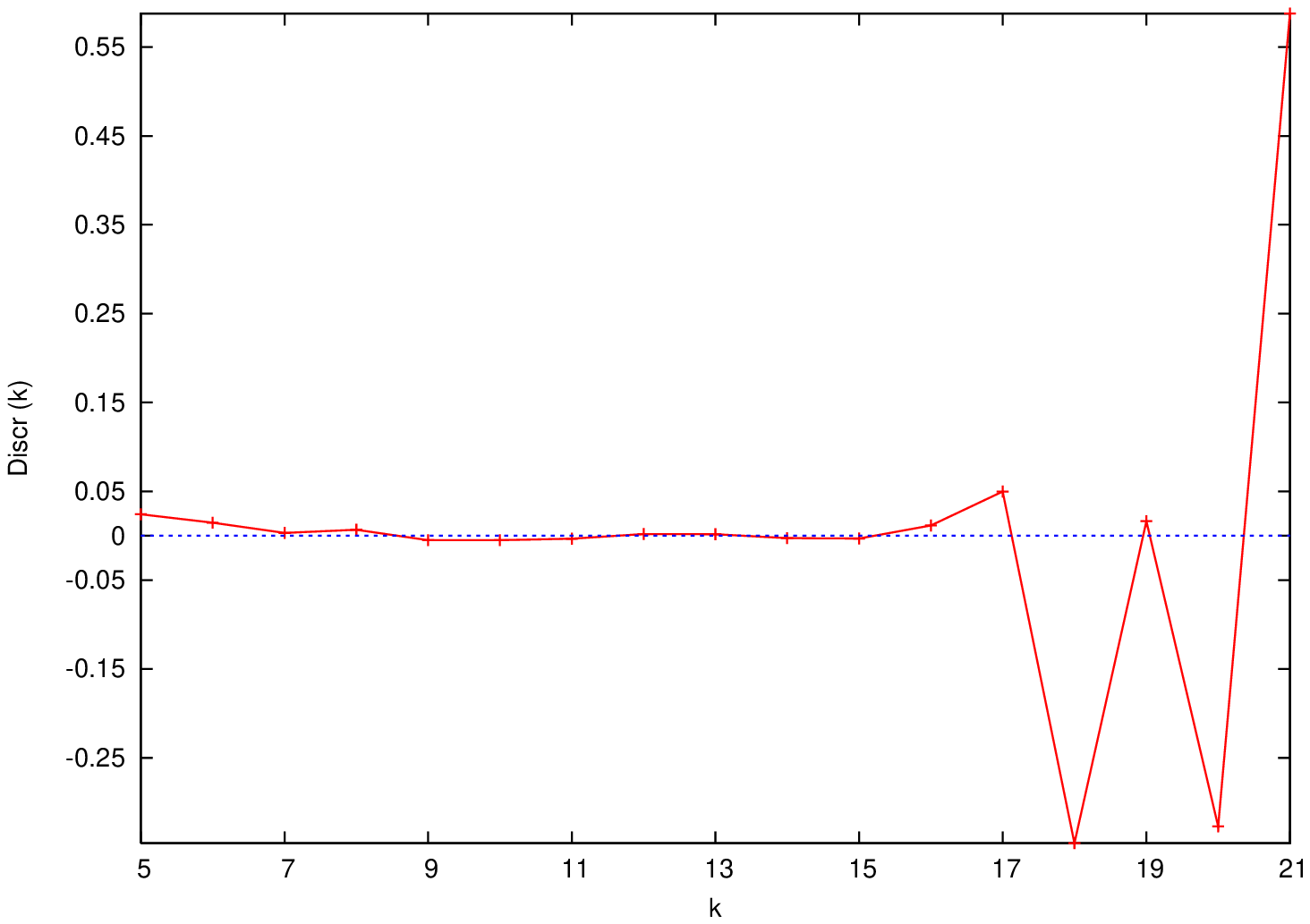}%
}
\caption{Discrepancy Discr vs wavenumber $k$ for the same data
as in Fig.~\ref{f:stupid-discrepancy-spectral}, but processed by a 5-stage
asymptotic extrapolation method.}
\label{f:joris-discrepancy-spectral}
\end{figure}
In the present case, the appropriate sequence of transformations is: {\bf
  Log}, {\bf D}, {\bf D}, {\bf I}, {\bf D}.  Because of the relatively
  low precision of the data it is not possible to perform more than
  five transformations, so that the method gives us access only to the
  leading-order asymptotic behavior, namely $|\hat u(k,1)| \simeq \ue
  ^{-C_\star|k|\ln |k|}$. It may be shown that the constant $C_\star=
  -1/u^{\rm (5)}$ where $u^{\rm (5)}$ is the constant value of the high-$|k|$
  interpolation $u^{\rm (5)}(|k|)$ after the 5th stage of transformation.
  Fig.~\ref{f:joris-discrepancy-spectral} shows the discrepancy $u^{\rm (5)}(|k|)+\ln 2$ . The absolute value
  of the 
  discrepancy lingers around 0.002 to 0.005 before being swamped by
  rounding noise at wavenumber 17.  Thus with asymptotic extrapolation
  the discrepancy does not exceed 0.7\% of the nominal value $\ln 2$.
  The accuracy of the determination has thus improved by about a
  factor 5, compared to the naive method without asymptotic
  extrapolation.

To improve further on this result and get some indication as to the
type of subdominant corrections present in the large-wavenumber
expansion of the Fourier coefficients it would not suffice to
increase the spatial resolution, since rounding noise would still
swamp the signal beyond a wavenumber of roughly 17. Higher precision
spectral calculations are doable but not very simple because high-precision
fast Fourier transform packages are still in the experimental phase.

In the next Section we shall present an alternative method,
significantly less versatile as to the choice of the initial condition
because it exploits the algebraic structure of a certain special class
of solutions, but which also allows to work easily in arbitrary
precision and thus to make better use of asymptotic extrapolation for
determining the constant $C_\star$.


\subsection{Half-space (Fourier) supported initial conditions}
\label{ss:simulation2}

%
%
%

So far we have limited ourselves to initial conditions that are real
entire functions. Hence the Fourier coefficients had Hermitian
symmetry: $\hat u_0 (-k) = \hat u_0^*(k)$, where the star denotes
complex conjugation.  With complex initial data there are no
analytical results when the dissipation is exponential, even when the
initial conditions are entire because the energy conservation relation
-- now about a complex-valued quantity -- ceases to give $L^2$-type
bounds. Actually, as already pointed out, Li and Sinai \cite{LS07}
showed that the 3D NSE can display finite-time blow up with suitable
complex initial data. It is however straightforward to adapt to
complex solutions the heuristic argument of
Section~\ref{ss:heuristics} and to predict a high-wavenumber
leading-order term exactly of the same form as for real
solutions. This is of interest since we shall see that there is a
class of periodic complex initial conditions for which, provided the
Burgers equation is written in terms as the Fourier coefficients as in
\cite{Sinai05} and \cite{LS07}, any given Fourier coefficient can be
calculated at arbitrary times $t$ with a finite number of operations,
most easily performed on a computer, by using either symbolic
manipulations or arbitrary-high precision floating point calculations.

For the case of the Burgers equation, this class consists of initial
conditions having the Fourier coefficients supported in the half line $k>0$.\footnote{If the coefficient for wavenumber $k=0$ is
non-vanishing a simple Galilean transformation can be used
to make it vanish.} We shall refer to this class of  
initial data as ``half-space (Fourier) supported''.

Because of the convolution structure of the nonlinearity when written
in terms of Fourier coefficients, it is obvious that with an
initial condition supported in the $k>0$ half line, the solution
will also be supported in this half line. A similar idea has been used 
in three dimensions for studying the singularities for complex
solutions of the 3D Euler equations  \cite{caflisch}.

Specifically, we consider again the 1D Burgers equation
\eqref{burgers}-\eqref{initburgers}
with $2\pi$-periodic boundary conditions  for $t\ge 0$, rewritten as
\eqref{fourierburgers},
in terms of the Fourier
coefficients $\hat u(k,t)$,  assumed here to exist.
The $k$-dependent real, even, non-negative
dissipation coefficient $\mu \ue ^{G(k)}$ is denoted by $\rho(k)$. The
initial conditions $\hat u_0(k)$ ($k=1,2,\ldots$) are chosen
arbitrarily, real or complex. We then have the following Proposition,
which is of purely algebraic nature:

\textbf{Proposition~5.1} \,\, {\it Eq.~\eqref{fourierburgers}
  with the initial conditions $\hat u(k,0) =0$ ($k\le 0$) and $\hat u(k,0)=\hat
  u_0(k)$ ($k=1,2,\ldots$) defines, for all $k>0$ and $t>0$,  $\hat
  u(k,t)$ as 
a polynomial functions of the set of 
$\hat u_0(p)$ with $0<p\le k$.
}

\textit{Proof}  \,\, From \eqref{fourierburgers},  after 
integration of the dissipative term, we obtain, for $k>0$ and $t>0$
\begin{equation}
\hat u(k,t) = \ue ^{-t \rho(k)} \hat u_0(k) -\frac{\ui k}{2}\int_0^t ds\,
\ue ^{-(t-s)\rho(k)}\sum_{p=1}^{k-1} \hat
u(p,s)\hat u(k-p,s).\label{neumann}
\end{equation}
Observe that $\hat u(k,t)$, given
by \eqref{neumann}, involves $\hat u(k,0)$ (linearly) and the  set of
Fourier
coefficients $\hat
u(p,s)$ for $1\le p\le k-1$ and $0\le s<t$ (quadratically). The proof
follows
by recursive use of this property for
$k$, $k-1$, \ldots, $1$.\footnote{This 
proposition has an obvious counterpart for the NSE in any dimension
when
the the Fourier coefficients of the initial condition are 
compactly supported in a product of half-spaces.} 

Note that the solution can be obtained without
any truncation error on a computer, using symbolic
manipulation. Alternatively it can be calculated in arbitrary high-precision
floating point arithmetic.

Now we specialize the initial condition even further, by assuming
that the only Fourier harmonic present in the initial condition has
$k=1$.\footnote{What follows can be easily extended to the case of a
  finite number of non-vanishing initial Fourier harmonics.}
We take
\begin{equation}
u_0(x) = \ui {\cal A} \ue ^{ \ui x},
\label{halfmode}
\end{equation}
for which $\hat u_0(1) = \ui {\cal A}$, while all the other coefficients
vanish. 

Setting, for $k>0$,
\begin{equation}
\hat v(k,t) \coloneqq \frac{\hat u(k,t)}{\ui  {\cal A} ^{k}},
\label{defv}
\end{equation}
we obtain from \eqref{neumann} by working out the power series to
the second order the following fully explicit expressions of the first
two  Fourier coefficient at any time $t>0$:
\begin{equation}
\hat{v} (1,t) =  \ue^{- \rho(1) t } , \qquad 
\hat{v} (2,t) =  \frac{\ue^{-2 \rho(1) t } }{\rho(2) - 2 \rho(1) }   -   \frac{\ue^{- \rho(2) t } }{\rho(2) - 2 \rho(1) }.  \label{firstthree}
\end{equation}

For the one-mode initial condition with ${\cal A}=1$ and $\rho(k) =  \ue
^{|k|}$ (exponential dissipation with $\mu =1$  and ($\kd =1$),
we have calculated the Fourier coefficients $\hat{v} (k,1) $ for $k=1,
2,\ldots,24$ 
using Maple symbolic calculation with a forty-digit accuracy.

\begin{figure}
\centerline{%
\includegraphics[scale=0.75]{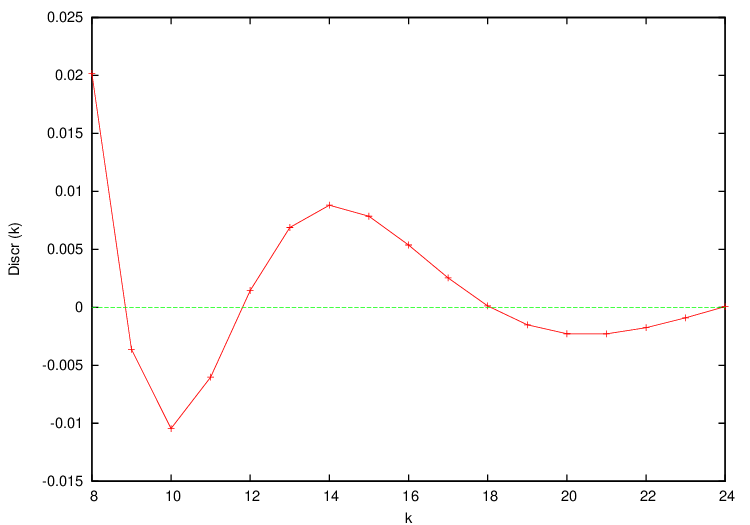}%
}
\caption{Fifth stage of asymptotic interpolation showing the 
discrepancy with respect to the heuristic prediction $C_\star =
  1/\ln2$. Note the decaying oscillatory behavior.} 
\label{f:joris-discrepancy-sinai}
\end{figure}

The data have then been processed using asymptotic extrapolation 
as in Section~\ref{ss:simulation1} with the same five transformations
{\bf Log}, {\bf D}, {\bf D}, {\bf I}, {\bf
  D}. ~Fig.~\ref{f:joris-discrepancy-sinai} shows the discrepancy 
$\hat v^{\rm (5)}(k,1)+\ln 2$ between
the 5th stage of interpolation and the prediction from the dominant 
balance argument. It is seen that the discrepancy drops to 
$-2.3\times 10^{-3}$. Thus the relative error is about $3\times 10^{-3}$.
The oscillations in the discrepancy, if they continue to higher 
wavenumbers would indicate that the first subdominant correction
to  the asymptotic behavior of the Fourier coefficient $\ue ^{-(1/\ln
  2) |\kt \ln |\kt|}$  is a prefactor involving a complex power of the
wavenumber.  

Finally, we address the issue of what kind of solution we have
 constructed by this Fourier-based algebraic method. Is it a
 ``classical'' sufficiently smooth global-in-time
solution of the Burgers equations \eqref{burgers}-\eqref{initburgers}
 written in the physical space? We have here obtained strong numerical
 evidence
that the Fourier coefficients decrease faster than exponentially with
the wavenumber and thus define a classical solution which is an
entire function of the space variable. This is however just a
 conjecture. The tools used in Section~\ref{s:entire} to prove the
 entire
character of the solution rely heavily on the definite positive character
 of the energy, a property lost with complex solutions.\footnote{It is 
however not difficult to prove, for short times, that our solution is
 also a  classical  entire solution.}



\section{Conclusions}
\label{s:conclusion}


In this paper we have proved that for a large class of evolution PDE's,
including the 3D NSE, exponential or faster-growing dissipation implies that
the solution becomes and remains an  entire function in the space variables at all
times. Exponential
growth constitutes a threshold: subexponential growth with  a Fourier
symbol
$\ue ^{|k|^\alpha}$, where $0<\alpha <1$ makes the solution analytic (but
not entire) as is the case in 2D (and generally conjectured in 3D).
Furthermore,
for the 3D NSE and the 1D Burgers equation with a dissipation having the
Fourier symbol $\mu \ue ^{|k|/\kd}$, we have shown that the amplitude of the
Fourier coefficients is bounded by $\ue ^{- C|\kt|\ln |\kt|}$ (where $\kt
\coloneqq k/\kd$) for any $C<1/(2\ln 2)$. For the case of the 1D Burgers
equation we have good evidence that this can be improved to $C<C_\star
= 1/\ln 2$  since the high-$|k|$ asymptotics seems to have a leading term precisely
of the form $\ue ^{- C_\star |\kt|\ln |\kt|}$; the evidence comes both from a
heuristic dominant balance argument and from high-precision simulations. The
heuristic argument can actually be carried over somewhat loosely to the expNSE
in any dimension: again the dominant nonlinear interaction contributing to
wave vector $k$ comes from the wave vectors $p = q =k/2$; actually the condition of
incompressibility kills nonlinear interactions between exactly parallel wave
vectors but this is only expected to modify algebraic prefactors in front of
the exponential term. 

We thus conjecture that the $C_\star = 1/\ln2$ also holds for expNSE
in any space dimension $d\ge 2$. Of course there is a substantial gap
between the bound and the conjectured asymptotic behavior. It seems
that such a gap is hard to avoid when using $L^2$-type norms. For
proving the entire character of the solution such norms were
appropriate. Beyond this, it is appears more advisable to try bounding
directly the moduli of Fourier coefficients by using the power series
method \cite{LS07,Sinai05}. A first step in this direction would be to
prove that $C_\star = 1/\ln2$ for initial conditions whose Fourier
coefficients are compactly supported in a product of half-spaces, of
the kind considered in Section~\ref{ss:simulation2}.\footnote{Progress on
this issue has been made and will be reported elsewhere.}

Exponential dissipation differs from ordinary dissipation
(with a Laplacian or a power thereof) not only by giving a
faster decay of the Fourier coefficients but by doing so
in a \emph{universal} way: with ordinary dissipation the decay
of the Fourier coefficient is generally conjectured to be, to leading
order,
of the form $\ue ^ {-\eta |k|}$ where $\eta$ depends on the viscosity
$\nu$ and on the energy input or on the size of the initial velocity;
with exponential dissipation the decay is $\ue ^{- C_\star |\kt|\ln
  |\kt|}$ where $C_\star =1/\ln 2$ and thus depends neither on the coefficient
$\mu$ which plays the role of the viscosity nor on the initial data.\footnote{It does however depend on the type of
  nonlinearity. For example with a cubic nonlinearity the same kind of
  heuristics as presented in Section~\ref{ss:heuristics} predicts a
  constant $C_\star =1/\ln 3$.}
As a consequence, it is expected that exponential dissipation 
will not exhibit the phenomenon of dissipation-range intermittency,
which for the usual dissipation can be traced back either to the fluctuations of $\eta$
\cite{rhk67} or to complex singularities of a velocity field that is
analytic  but not entire \cite{frischmorf}. 

Finally some comments on the practical relevance of modified
dissipation.   First, let us comment on ``hyperviscosity'', the replacement
of the (negative) Laplacian by its power of order $\alpha>1$.
Of course we know that 
specialists of PDE's have
traditionally been interested in the hyperviscous 3D NSE, 
perhaps to overcome the frustration of not being able
to prove much about the ordinary 3D NSE. But scientists doing numerical
simulations of the NSE, say, for engineering, astrophysical or
geophysical applications, have also been using hyperviscosity 
because it is often believed to  allow effectively higher
Reynolds numbers without the need to increase spatial resolution. 
Recently, three of us (UF, WP, SSR) and other coauthors have shown
that when using a high power $\alpha$ of the Laplacian in the dissipative 
term for 3D NSE or 1D Burgers, one comes very close to a Galerkin
truncation of Euler or inviscid Burgers, respectively \cite{FKPPRWZ08}. This produces
a range of nearly thermalized modes which shows up in  large-Reynolds
number spectral simulations
as a huge bottleneck in the Fourier amplitudes between the inertial
range and  the far dissipation range. Since the bottleneck
generates a fairly large eddy viscosity, the hyperviscosity procedure
with large $\alpha$ actually \emph{decreases} the effective Reynolds 
number. 

Next, consider exponential dissipation. In 1996 Achim Wirth noticed that when
used in the 1D Burgers equation, cosh dissipation produces almost no
bottleneck although it grows much faster than a power of the wavenumber at high
wavenumbers \cite{AW}.  It is now clear that such a dissipation will
produce
a faster-than-exponential decay at the highest wavenumbers. But at 
wavenumbers such that $|k|\ll \kd$ a dissipation rate $-\mu (1-\cosh k/\kd)$
reduces to $\mu|k|^2/(2\kd ^2)$, to leading order, which is the
ordinary (Laplacian) dissipation. With 
 the ordinary 1D Burgers equation it may be shown analytically that
there is no bottleneck.  For the ordinary 3D NSE, experimental and numerical results
show the presence of a rather modest bottleneck (for example the
``compensated'' three-dimensional energy spectrum $|k|^{+5/3} E(|k|)$ overshoots by about
20\%.).  If in a simulation with cosh dissipation $\mu$ and $\kd$ are adjusted
in such a way that dissipation starts acting at wavenumbers slighter smaller
than $\kd$, the beginning of the dissipation range will be mostly as with an
ordinary Laplacian, that is with no or little bottleneck.\footnote{If $\mu$
and $\kd$ are not carefully chosen, effective dissipation can start well
beyond $\kd$.  One may then observe the same kind of thermalization and of
bottleneck than with a high power of the Laplacian \cite{JZZ}.}  At higher
wavenumbers, where the exponential growth of the dissipation rate is felt,
faster than exponential decay will be observed.  In principle this can be used
to avoid wasting resolution without developing a serious bottleneck.  Faster
than exponentially growing dissipation, e.g. $\mu \left(\ue ^{(|k|/\kd)^2}
-1\right)$, may be even better because the prediction is that the Fourier
coefficients will display Gaussian decay.\footnote{Here we mention
  that this  may be of relevance for a numerical procedure where a Gaussian
  filter is used   at
each time step, a procedure described to one of us (UF) as allowing
to absorb energy near the maximum wavenumber without having it
reflected back to lower wavenumbers \cite{orszag}.}   Testing the advantages and
drawbacks of different types of faster-than-algebraically growing dissipations
for numerical simulations is left for future work.


%
%
%

We thank J.-Z.~Zhu and A.~Wirth for important input and M.~Blank, K.~Khanin,
B.~Khesin and V.~Zheligovsky for many remarks.
CB acknowledges the warm hospitality of the Weizmann
Institute and SSR that of the Observatoire de la C\^ote d'Azur, places where 
parts of this work were carried out. The work of EST was supported in part
by
the NSF grant No.~DMS-0708832 and the ISF grant No.~120/06 .
SSR  thanks R. ~Pandit, D.~Mitra and P.~Perlekar for useful discussions
and acknowledges DST and UGC (India) for
support and SERC (IISc) for computational resources. UF, WP and SSR 
were partially supported by ANR ``OTARIE'' BLAN07-2\_183172
and used  the M\'esocentre de calcul of the Observatoire de la C\^ote
d'Azur for computations.

\end{document}